\documentclass[10pt]{article}
\usepackage{amsfonts}
\usepackage{amsmath}
\usepackage{amscd,amssymb,amsthm}
\usepackage{graphics}
\usepackage{graphicx}

\newtheorem{theorem}{Theorem}[section]

\newtheorem{example}[theorem]{Example}

\newtheorem{proposition}[theorem]{Proposition}
\newtheorem{remark}[theorem]{Remark}

\begin{document}

\title{Degenerate Chenciner bifurcation revisited}
\author{G. Moza\thanks{
Department of Mathematics, Politehnica University of Timisoara, Romania;
email: gheorghe.moza@upt.ro}, O. Brandibur\thanks{
Department of Mathematics, West University of Timisoara, Romania}, E.A.
Kokovics, L.F. Vesa}

\date{}
\maketitle

\begin{abstract}
Generic results for degenerate Chenciner (generalized Neimark-Sacker)
bifurcation are obtained in the present work. The bifurcation arises in
two-dimensional discrete-time systems with two independent parameters. We
define in this work a new transformation of parameters, which enables the
study of the bifurcation when the degeneracy occurs. By the four bifurcation
diagrams we obtain, new behaviors hidden by the degeneracy are brought to
light.
\end{abstract}


\section{Introduction}

\quad\ Many real-world applications are modeled using dynamical systems
described by both differential equations and their discrete-time
counterparts, namely, difference equations. A discrete-time approach for
studying various models is proved to be efficient, especially from a
computational point of view. They may uncover complex behaviors which are
not easily captured by a continuous-time approach. A review on
continuous-time versus discrete-time approaches in scheduling of chemical
processes is presented in \cite{flo}. The study of discrete-time dynamical
systems is an active field of research \cite{bro1}, \cite{det1}--\cite{kha1}%
, \cite{lib1}--\cite{vit1}. The analysis of bifurcations is among the most
relevant topics in the qualitative theory of dynamical systems,
particularly, in discrete-time systems \cite{wig}.

In this work, we focus our attention on the Chenciner bifurcation (known
also as the generalized Neimark-Sacker bifurcation) [Chenciner, 1985, 1988], 
\cite{che3}, which arises in two-dimensional discrete-time dynamical systems
with two real parameters. The presence of this bifurcation in a three-mass
chain model has been reported in \cite{bakri}. The non-degenerate case of
this bifurcation has been studied earlier \cite{Ku98}, while a degenerate
case has been recently presented in \cite{tig2}.

The purpose of this paper is to further study the degenerate case from \cite%
{tig2}, by considering a different approach, which reduces the number of
bifurcation diagrams that are needed to describe the bifurcation in the
degenerate case.

The paper is organized as it follows. After the introduction, in section two
we present the main ingredients of the bifurcation and establish the
research objectives to be explored later in the work. Section three presents
briefly the general results which are needed for studying the bifurcation in
a degenerate framework. The main results of the paper are comprised in
section four, where we define a new transformation of parameters for the
degeneracy we study. Using this change we obtain four bifurcation diagrams,
which describe the behavior of the system when the degeneracy occurs. The
last section contains conclusive remarks on our results.

\section{Preliminaries}

Consider the discrete-time system 
\begin{equation}
x_{n+1}=f\left( x_{n},\alpha \right) ,  \label{ecn}
\end{equation}%
where $x_{n}\in 
\mathbb{R}
^{2},$ $n\in 
\mathbb{Z}
,$ $\alpha =\left( \alpha _{1},\alpha _{2}\right) \in 
\mathbb{R}
^{2}$ and $f$ is a smooth function of class $C^{r}$ with $r\geq 2.$ Another
form of (\ref{ecn}) which avoids indices is 
\begin{equation}
x\longmapsto f\left( x,\alpha \right) .  \label{ecs}
\end{equation}%
Using complex coordinates, (\ref{ecs}) can be written as 
\begin{equation}
z\longmapsto \delta \left( \alpha \right) z+g\left( z,\bar{z},\alpha \right)
,  \label{ecz}
\end{equation}%
where $\delta $ and $g$ are smooth functions, $\delta
\left( \alpha \right) =r\left( \alpha \right) e^{i\theta \left( \alpha
\right) }$ and $g\left( z,\bar{z},\alpha \right) =\sum_{k+l\geq 2}\frac{1}{k!l!}%
g_{kl}\left( \alpha \right) z^{k}\bar{z}^{l},$ with $r\left( 0\right) =1$
and $\theta \left( 0\right) =\theta _{0};$ $g_{kl}\left( \alpha \right) $
are smooth complex-valued functions.

Equation (\ref{ecz}) can be further transformed as in \cite{Ku98} to 
\begin{equation}
w\longmapsto \left( r\left( \alpha \right) +b_{1}\left( \alpha \right) w\bar{%
w}+b_{2}\left( \alpha \right) w^{2}\bar{w}^{2}\right) we^{i\theta \left(
\alpha \right) }+O\left( \left\vert w\right\vert ^{6}\right) ,  \label{ecw}
\end{equation}%
where $b_{k}\left( \alpha \right) =a_{k}\left( \alpha \right) e^{-i\theta
\left( \alpha \right) },$ $k=1,2.$

\begin{remark}
Equation (\ref{ecz}) becames (\ref{ecw}) by using the following invertible
smoothly parameter change of complex coordinate 
\begin{equation*}
z=w+\sum\limits_{2\leq k+l\leq 5}^{{}}\dfrac{1}{k!l!}h_{kl}(\alpha )w^{k}%
\bar{w}^{l},\quad h_{21}(\alpha )=h_{32}(\alpha )=0.
\end{equation*}
\end{remark}

Denoting by 
\begin{equation}
\beta _{1}\left( \alpha \right) =r\left( \alpha \right) -1\text{ and }\beta
_{2}\left( \alpha \right) =Re\left( b_{1}(\alpha )\right) ,  \label{b12}
\end{equation}%
and using polar coordinates, (\ref{ecw}) becomes

\begin{equation}
\left\{ 
\begin{array}{cc}
\rho _{n+1}= & \rho _{n}\left( 1+\beta _{1}\left( \alpha \right) +\beta
_{2}\left( \alpha \right) \rho _{n}^{2}+\allowbreak L_{2}\left( \alpha
\right) \rho _{n}^{4}\right) +\rho _{n}O\left( \rho _{n}^{6}\right) \\ 
\varphi _{n+1}= & \varphi _{n}+\theta \left( \alpha \right) +\rho
_{n}^{2}\left( \frac{Im\left( b_{1}\left( \alpha \right) \right) }{\beta
_{1}\left( \alpha \right) +1}+O\left( \rho _{n},\alpha \right) \allowbreak
\right) \ \ \ \ \ \ \ \ \ \text{ }%
\end{array}%
\right. ,  \label{rofi}
\end{equation}%
where $L_{2}\left( \alpha \right) =\frac{Im^{2}\left( b_{1}\left( \alpha
\right) \right) +2\left( 1+\beta _{1}\left( \alpha \right) \right) Re\left(
b_{2}\left( \alpha \right) \right) }{2\left( \beta _{1}\left( \alpha \right)
+1\right) }.$

A bifurcation in the system (\ref{rofi}) satisfying $r\left( 0\right) =1,$ $%
Re\left( b_{1}(0)\right) =0$ and $L_{2}\left( 0\right) \neq 0,$ is known as
the \textit{Chenciner bifurcation}. Since $\beta _{1}\left( 0\right) =0,$ it
follows that 
\begin{equation*}
L_{2}\left( 0\right) =\frac{1}{2}\left( Im^{2}\left( b_{1}\left( 0\right)
\right) +2Re\left( b_{2}\left( 0\right) \right) \right) .
\end{equation*}

When the transformation of parameters 
\begin{equation}
\left( \alpha _{1},\alpha _{2}\right) \longmapsto \left( \beta _{1}\left(
\alpha \right) ,\beta _{2}\left( \alpha \right) \right)  \label{trp}
\end{equation}%
is regular at $\left( 0,0\right) ,$ then $\beta _{1}$ and $\beta _{2}$
become the new parameters of the system (\ref{rofi}). This is the \textit{%
non-degenerate case} of the bifurcation.

The aim of this paper is to study the bifurcation in the \textit{degenerate}
case with respect to the transformation of parameters (\ref{trp}), namely,
when (\ref{trp}) is not regular at $\left( 0,0\right) .$ Hence, we cannot
use $\beta _{1,2}$ as new parameters of the system (\ref{rofi}), due to this
degeneracy condition. Alternatively, we could keep on working with the
initial parameters $\alpha _{1,2}$ in the polar form (\ref{rofi}), as it is
done in \cite{tig2}, or explore for another regular transformation of the
parameters. We opt for the second method in this paper.

\section{Analysis of degenerate Chenciner bifurcation}

The truncated form of the $\rho -$map of (\ref{rofi}) obtained by
eliminating higher order terms is 
\begin{equation}
\rho _{n+1}=\rho _{n}\left( 1+\beta _{1}\left( \alpha \right) +\beta
_{2}\left( \alpha \right) \rho _{n}^{2}+\allowbreak L_{2}\left( \alpha
\right) \rho _{n}^{4}\right) .  \label{tru}
\end{equation}

The $\varphi -$map of (\ref{rofi}), which describes a rotation by an angle
depending on $\alpha $ and $\rho ,$ can be approximated by its truncated
form 
\begin{equation}
\varphi _{n+1}=\varphi _{n}+\theta \left( \alpha \right) .  \label{truf}
\end{equation}

Assume $0<\theta \left( 0\right) <\pi .$ Henceforward, the system to be
studied in this paper is (\ref{tru})-(\ref{truf}), which is known as \textit{%
the (truncated) normal form} of the system (\ref{ecw}). The main equation
shaping the dynamics of the system (\ref{tru})-(\ref{truf}) is the $\rho -$%
map (\ref{tru}), which is independent from the $\varphi -$map and therefore,
will be studied separately.

Defining a one-dimensional dynamical system, the $\rho -$map has the fixed
point $\rho =0,$ for all values of $\alpha ,$ which corresponds to the fixed
point $O\left( 0,0\right) $ in the normal form (\ref{tru})-(\ref{truf}). A
positive nonzero fixed point of the $\rho -$map corresponds to an invariant
closed curve in (\ref{tru})-(\ref{truf}).

Observe that $sign\left( L_{2}\left( \alpha \right) \right) =sign\left(
L_{0}\right) $ for $\left\vert \alpha \right\vert =\sqrt{\alpha
_{1}^{2}+\alpha _{2}^{2}}$ sufficiently small, as 
\begin{equation*}
L_{2}\left( \alpha \right) =L_{0}\left( 1+O\left( \left\vert \alpha
\right\vert \right) \right) \quad \text{and}\quad L_{0}\neq 0.
\end{equation*}

Throughout the work, $O\left( \left\vert \alpha \right\vert ^{n}\right) $
denotes the higher order terms in a Taylor expansion at $\alpha =0.$ For
example, 
\begin{equation*}
O\left( \left\vert \alpha \right\vert \right) =k_{10}\alpha
_{1}+k_{01}\alpha _{2}+...\ \ \ \ \text{ and }\ \ \ \ O\left( \left\vert
\alpha \right\vert ^{2}\right) =k_{20}\alpha _{1}^{2}+k_{11}\alpha
_{1}\alpha _{2}+k_{02}\alpha _{2}^{2}+...
\end{equation*}

The next Proposition \ref{pr1} describes the stability of $O$ for $%
\left\vert \alpha \right\vert $ sufficiently small, while Theorem \ref{th1}
deals with the existence of invariant closed curves in the normal form (\ref%
{tru})-(\ref{truf}). Their proofs can be obtained by studying (\ref{tru})
and are presented in \cite{tig2}.

\begin{proposition}
\label{pr1} The fixed point $O$ is (linearly) stable if $\beta _{1}\left(
\alpha \right) <0$ and unstable if $\beta _{1}\left( \alpha \right) >0,$ for
all values $\alpha $ with $\left\vert \alpha \right\vert $ sufficiently
small. On the bifurcation curve $\beta _{1}\left( \alpha \right) =0,$ $O$ is
(nonlinearly) stable if $\beta _{2}\left( \alpha \right) <0$ and unstable if 
$\beta _{2}\left( \alpha \right) >0,$ when $\left\vert \alpha \right\vert $
is sufficiently small. At $\alpha =0,$ $O$ is (nonlinearly) stable if $%
L_{0}<0$ and unstable if $L_{0}>0.$
\end{proposition}

The positive nonzero fixed points of the $\rho -$map (\ref{tru}), which are
the solutions of the equation 
\begin{equation}
L_{2}\left( \alpha \right) y^{2}+\beta _{2}\left( \alpha \right) y+\beta
_{1}\left( \alpha \right) =0,  \label{y1y2}
\end{equation}%
where $y=\rho _{n}^{2},$ generate invariant closed curves (invariant
circles) in the system (\ref{tru})-(\ref{truf}).

Denote by 
\begin{equation}
\Delta \left( \alpha \right) =\beta _{2}^{2}\left( \alpha \right) -4\beta
_{1}\left( \alpha \right) L_{2}\left( \alpha \right) ,  \label{delta}
\end{equation}%
respectively, $y_{1}=\frac{1}{2L_{2}\left( \alpha \right) }\left( \sqrt{%
\Delta \left( \alpha \right) }-\beta _{2}\left( \alpha \right) \right) $ and 
$y_{2}=-\frac{1}{2L_{2}\left( \alpha \right) }\left( \sqrt{\Delta \left(
\alpha \right) }+\beta _{2}\left( \alpha \right) \right) $ the two roots of (%
\ref{y1y2}), whenever they exist as real numbers.

\begin{theorem}
\label{th1} The following assertions are true.

1) When $\Delta \left( \alpha \right) <0$ for all $\left\vert \alpha
\right\vert $ sufficiently small, the normal form (\ref{tru})-(\ref{truf})
has no invariant circles.

2) When $\Delta \left( \alpha \right) >0$ for all $\left\vert \alpha
\right\vert $ sufficiently small, the normal form (\ref{tru})-(\ref{truf})
has:

\qquad a) one invariant unstable circle $\rho _{n}=\sqrt{y_{1}}$ if $L_{0}>0$
and $\beta _{1}\left( \alpha \right) <0;$

\qquad b) one invariant stable circle $\rho _{n}=\sqrt{y_{2}}$ if $L_{0}<0$
and $\beta _{1}\left( \alpha \right) >0;$

\qquad c) two invariant circles, $\rho _{n}=\sqrt{y_{1}}$ unstable and $\rho
_{n}=\sqrt{y_{2}}$ stable, if $L_{0}>0,$ $\beta _{1}\left( \alpha \right)
>0, $ $\beta _{2}\left( \alpha \right) <0$ or $L_{0}<0,$ $\beta _{1}\left(
\alpha \right) <0,$ $\beta _{2}\left( \alpha \right) >0;$ in addition, $%
y_{1}<y_{2}$ if $L_{0}<0$ and $y_{2}<y_{1}$ if $L_{0}>0;$

\qquad d) no invariant circles if $L_{0}>0,$ $\beta _{1}\left( \alpha
\right) >0,$ $\beta _{2}\left( \alpha \right) >0$ or $L_{0}<0,$ $\beta
_{1}\left( \alpha \right) <0,$ $\beta _{2}\left( \alpha \right) <0.$

3) On the bifurcation curve $\Delta \left( \alpha \right) =0,$ the system (%
\ref{tru})-(\ref{truf}) has one invariant unstable circle $\rho _{n}=\sqrt{%
y_{1}}$ for all $L_{0}\neq 0.$

4) When $\beta _{1}\left( \alpha \right) =0,$ the system (\ref{tru})-(\ref%
{truf}) has one invariant circle $\rho _{n}=\sqrt{-\frac{\beta _{2}\left(
\alpha \right) }{L_{0}}}$ whenever $L_{0}\beta _{2}\left( \alpha \right) <0.$
It is stable if $L_{0}<0$ and $\beta _{2}\left( \alpha \right) >0,$
respectively, unstable if $L_{0}>0$ and $\beta _{2}\left( \alpha \right) <0.$
\end{theorem}

\section{Bifurcation diagrams}

The transformation (\ref{trp}) is not regular at $\alpha =0$ and, thus, the
Chenciner bifurcation is \textit{degenerate}, if and only if $\left. \frac{%
\partial \beta _{1}}{\partial \alpha _{1}}\frac{\partial \beta _{2}}{%
\partial \alpha _{2}}-\frac{\partial \beta _{1}}{\partial \alpha _{2}}\frac{%
\partial \beta _{2}}{\partial \alpha _{1}}\right\vert _{\alpha =0}=0.$
Denoting the coefficients of the linear parts of $\beta _{1,2}\left( \alpha
\right) $ by $\frac{\partial \beta _{1}}{\partial \alpha _{i}}\left(
0\right) =c_{i},$ respectively, $\frac{\partial \beta _{2}}{\partial \alpha
_{i}}\left( 0\right) =d_{i},$ $i=1,2,$ the degeneracy condition becomes%
\begin{equation}
c_{1}d_{2}-c_{2}d_{1}=0.  \label{dc1}
\end{equation}%
Assume $c_{i}\neq 0\ $and $d_{i}\neq 0.$ Denote by $L_{2}\left( \alpha
\right) =L_{0}+l_{1}\alpha _{1}+l_{2}\alpha _{2}+O\left( \left\vert \alpha
\right\vert ^{2}\right) .$ In general, denote by $\beta _{1}\left( \alpha
\right) =c_{1}\alpha _{1}+c_{2}\alpha
_{2}+\sum\limits_{i+j=2}^{m}c_{ij}\alpha _{1}^{i}\alpha _{2}^{j}+O\left(
\left\vert \alpha \right\vert ^{m+1}\right) $ and $\beta _{2}\left( \alpha
\right) =d_{1}\alpha _{1}+d_{2}\alpha
_{2}+\sum\limits_{i+j=2}^{n}d_{ij}\alpha _{1}^{i}\alpha _{2}^{j}+O\left(
\left\vert \alpha \right\vert ^{n+1}\right) .$

\begin{theorem}\label{th2} 
Assume the degeneracy condition (\ref{dc1}) holds true. Then, there
exists a transformation of parameters $\left( \alpha _{1},\alpha _{2}\right) 
\overset{S}{\longmapsto }\left( \mu _{1},\mu _{2}\right) $ given by 
\begin{equation}
\mu _{1}=\beta _{2}^{2}\left( \alpha \right) -4\beta _{1}\left( \alpha
\right) L_{2}\left( \alpha \right) \text{ and }\mu _{2}=\beta _{2}\left(
\alpha \right) +L_{2}\left( \alpha \right) -L_{0},  \label{tr2}
\end{equation}%
which is regular at $\alpha =0$ iff 
\begin{equation}
c_{1}l_{2}-c_{2}l_{1}\neq 0.  \label{newt}
\end{equation}%
Denote by $\widehat{\beta }_{1,2}\left( \mu \right) =\beta _{1,2}\circ
S^{-1}\left( \mu \right) .$ In the new parametric plane $\mu _{1}\mu _{2}$
and for $\left\vert \mu \right\vert $ sufficiently small, $\widehat{\beta }%
_{1}\left( \mu \right) =0$ is a curve of the form 
\begin{equation}
B_{1}=\left\{ \left( \mu _{1},\mu _{2}\right) \in 
\mathbb{R}
^{2},\mu _{1}=m_{2}^{2}\mu _{2}^{4}\left( 1+O\left( \mu _{2}\right) \right)
\right\} ,  \label{b1}
\end{equation}%
while $\widehat{\beta }_{2}\left( \mu \right) =0$ a curve given by 
\begin{equation}
B_{2}=\left\{ \left( \mu _{1},\mu _{2}\right) \in 
\mathbb{R}
^{2},\mu _{1}=4L_{0}\frac{c_{1}}{d_{1}}m_{2}\mu _{2}^{2}\left( 1+O\left( \mu
_{2}\right) \right) \right\} ,  \label{b2}
\end{equation}%
where $m_{2}$ is a real constant depending on the linear and quadratic
coefficients of $\beta _{1,2}\left( \alpha \right) $ and $L_{2}\left( \alpha
\right) .$ 
\end{theorem}

\textit{Proof.} For the transformation $S,$ denote by $A\left( \alpha
\right) =\left( 
\begin{array}{cc}
\frac{\partial \mu _{1}}{\partial \alpha _{1}}\left( \alpha \right) & \frac{%
\partial \mu _{1}}{\partial \alpha _{2}}\left( \alpha \right) \\ 
\frac{\partial \mu _{2}}{\partial \alpha _{1}}\left( \alpha \right) & \frac{%
\partial \mu _{2}}{\partial \alpha _{2}}\left( \alpha \right)%
\end{array}%
\right) .$ Using the linear terms of $\beta _{1,2}\left( \alpha \right) $
and $L_{2}\left( \alpha \right) ,$ the matrix $A\left( \alpha \right) $ at $%
\alpha =0$ becomes $A_{0}=\left( 
\begin{array}{cc}
-4L_{0}c_{1} & -4L_{0}c_{2} \\ 
d_{1}+l_{1} & d_{2}+l_{2}%
\end{array}%
\right) .$

The inverse transformation $\left( \mu _{1},\mu _{2}\right) \overset{S^{-1}}{%
\longmapsto }\left( \alpha _{1},\alpha _{2}\right) $ in its linear terms can
be determined from $\left( \mu _{1},\mu _{2}\right) ^{T}=A_{0}\cdot \left(
\alpha _{1},\alpha _{2}\right) ^{T},$ that is, $\left( \alpha _{1},\alpha
_{2}\right) ^{T}=A_{0}^{-1}\left( \mu _{1},\mu _{2}\right) ^{T}.$ We obtain 
\begin{equation}
\alpha _{1}=s_{10}\mu _{1}+s_{01}\mu _{2}\text{ and }\alpha _{2}=p_{10}\mu
_{1}+p_{01}\mu _{2},  \label{a12}
\end{equation}%
where $s_{10}=-\frac{d_{2}+l_{2}}{4L_{0}n_{0}},$ $s_{01}=-\frac{c_{2}}{n_{0}}%
,$ $p_{10}=\frac{d_{1}+l_{1}}{4L_{0}n_{0}},$ $p_{01}=\frac{c_{1}}{n_{0}},$
and $n_{0}=c_{1}d_{2}-c_{2}d_{1}+c_{1}l_{2}-c_{2}l_{1}.$

The inverse transformation $S^{-1}$ can be determined in further order
terms. For example, $S^{-1}$ in quadratic terms is of the form 
\begin{equation}
\alpha _{1}=\sum\limits_{i+j=1}^{2}s_{ij}\mu _{1}^{i}\mu _{2}^{j}\text{ and }%
\alpha _{2}=\sum\limits_{i+j=1}^{2}p_{ij}\mu _{1}^{i}\mu _{2}^{j}.
\label{q2}
\end{equation}%
The coefficients $s_{ij}$ and $p_{ij}$ can be determined in terms of the
coefficients of $\beta _{1,2}\left( \alpha \right) $ and $L_{2}\left( \alpha
\right) $ by the method of undetermined coefficients in (\ref{tr2}).

It follows from (\ref{dc1}) that 
\begin{equation*}
\det \left( A_{0}\right) =-4L_{0}\left( c_{1}l_{2}-c_{2}l_{1}\right) ,
\end{equation*}%
thus, since $L_{0}\neq 0,$ the new transformation $S$ is regular at $\alpha
=0$ if and only if $c_{1}l_{2}-c_{2}l_{1}\neq 0.$

Using (\ref{dc1}) in (\ref{a12}), $L_{2}\left( \alpha \right)
=L_{0}+l_{1}\alpha _{1}+l_{2}\alpha _{2}+O\left( \left\vert \alpha
\right\vert ^{2}\right) $ becomes $\widehat{L}_{2}\left( \mu \right)
=L_{2}\circ S^{-1}\left( \mu \right) $ given by 
\begin{equation}
\widehat{L}_{2}\left( \mu \right) =L_{0}+\frac{d_{1}}{4L_{0}c_{1}}\mu
_{1}+\allowbreak \mu _{2}+O\left( \left\vert \mu \right\vert ^{2}\right) ,
\label{el2}
\end{equation}%
respectively, $\widehat{\beta }_{2}\left( \mu \right) =\beta _{2}\circ
S^{-1}\left( \mu \right) $ by (\ref{tr2}) reads 
\begin{equation*}
\widehat{\beta }_{2}\left( \mu \right) =\mu _{2}-\widehat{L}_{2}\left( \mu
\right) +L_{0}=-\frac{1}{4}\frac{d_{1}}{L_{0}c_{1}}\mu _{1}+O\left(
\left\vert \mu \right\vert ^{2}\right) .
\end{equation*}%
In order to determine $\widehat{\beta }_{2}\left( \mu \right) $ in its
lowest terms, write $\widehat{L}_{2}\left( \mu \right) $ from (\ref{el2}) in
the form 
\begin{equation*}
\widehat{L}_{2}\left( \mu \right) =L_{0}+\frac{d_{1}}{4L_{0}c_{1}}\mu
_{1}\left( 1+O\left( \left\vert \mu \right\vert \right) \right) +\allowbreak
\mu _{2}-m_{2}\mu _{2}^{2}\left( 1+O\left( \mu _{2}\right) \right) ,
\end{equation*}%
where $m_{2}$ is assumed nonzero, $m_{2}\neq 0.$ Then, 
\begin{equation}
\widehat{\beta }_{2}\left( \mu \right) =-\frac{d_{1}}{4L_{0}c_{1}}\mu
_{1}\left( 1+O\left( \left\vert \mu \right\vert \right) \right) +m_{2}\mu
_{2}^{2}\left( 1+O\left( \mu _{2}\right) \right) .  \label{be2}
\end{equation}%
Since $\widehat{\beta }_{2}\left( 0,0\right) =0$ and $\frac{\partial 
\widehat{\beta }_{2}}{\partial \mu _{1}}\left( 0,0\right) =-\frac{1}{4}\frac{%
d_{1}}{L_{0}c_{1}}\neq 0,$ we can apply the Implicit Function Theorem to the
equation $\widehat{\beta }_{2}\left( \mu _{1},\mu _{2}\right) =0$ given by (%
\ref{be2}). Thus, there exists $\varepsilon >0$ sufficiently small and a
function 
\begin{equation*}
\mu _{1}:\left( -\varepsilon ,\varepsilon \right) \rightarrow 
\mathbb{R}
,\text{ }\mu _{1}=\mu _{1}\left( \mu _{2}\right) ,
\end{equation*}%
such that $\mu _{1}\left( 0\right) =0$ and $\widehat{\beta }_{2}\left( \mu
_{1}\left( \mu _{2}\right) ,\mu _{2}\right) =0$ for all $\left\vert \mu
_{2}\right\vert <\varepsilon ;$ $\mu _{1}$ becomes a function of argument $%
\mu _{2}$ in (\ref{be2}), $\mu _{1}=\mu _{1}\left( \mu _{2}\right) .$
Deriving now in $-\frac{d_{1}}{4L_{0}c_{1}}\mu _{1}\left( 1+O\left(
\left\vert \mu \right\vert \right) \right) +m_{2}\mu _{2}^{2}\left(
1+O\left( \mu _{2}\right) \right) =0$ with respect to $\mu _{2},$ where $\mu
_{1}=\mu _{1}\left( \mu _{2}\right) ,$ we obtain $\frac{\partial \mu _{1}}{%
\partial \mu _{2}}\left( 0\right) =0$ and $\frac{\partial ^{2}\mu _{1}}{%
\partial \mu _{2}^{2}}\left( 0\right) =8L_{0}\frac{c_{1}}{d_{1}}m_{2}.$

Expressing $\mu _{1}=\mu _{1}\left( 0\right) +\frac{\partial \mu _{1}}{%
\partial \mu _{2}}\left( 0\right) \mu _{2}+\frac{1}{2}\frac{\partial ^{2}\mu
_{1}}{\partial \mu _{2}^{2}}\left( 0\right) \mu _{2}^{2}+O\left( \mu
_{2}^{3}\right) $ as a Taylor series at $\mu _{2}=0,$ we obtain 
\begin{equation}
\mu _{1}=4L_{0}\frac{c_{1}}{d_{1}}m_{2}\mu _{2}^{2}\left( 1+O\left( \mu
_{2}\right) \right) .  \label{mu2}
\end{equation}

Similarly, by (\ref{dc1}), (\ref{a12}), (\ref{tr2}) and (\ref{be2}), $%
\widehat{\beta }_{1}\left( \mu \right) =-\frac{1}{4\widehat{L}_{2}\left( \mu
\right) }\left( \mu _{1}-\widehat{\beta }_{2}^{2}\left( \mu \right) \right) $
becomes

\begin{equation}
\widehat{\beta }_{1}\left( \mu \right) =-\frac{1}{4L_{0}}\left[ \mu
_{1}\left( 1+O\left( \left\vert \mu \right\vert \right) \right)
-m_{2}^{2}\mu _{2}^{4}\left( 1+O\left( \mu _{2}\right) \right) \right] ,
\label{be1}
\end{equation}%
which, by the Implicit Function Theorem, leads to $B_{1}.$

The exact expression of $m_{2}$ is not essential in the qualitative analysis
of the Chenciner bifurcation we aim to obtain in this article. However,
since in concrete applications it is useful, we determine it. To this end,
using (\ref{q2}) in $L_{2}\left( \alpha \right) =L_{0}+l_{1}\alpha
_{1}+l_{2}\alpha _{2}+l_{20}\alpha _{1}^{2}+l_{11}\alpha _{1}\alpha
_{2}+l_{02}\alpha _{2}^{2},\allowbreak $ it follows that the term $-m_{2}\mu
_{2}^{2}$ of $\widehat{L}_{2}\left( \mu \right) $ has the coefficient $%
m_{2}=-\left(
l_{02}p_{01}^{2}+l_{11}p_{01}s_{01}+l_{20}s_{01}^{2}+l_{2}p_{02}+l_{1}s_{02}%
\right) .$

We need further the inverse transformation (\ref{q2}), more exactly, we need
the coefficients $s_{ij}$ and $p_{ij}$ with $i+j=1$ (which are already
known), respectively, $s_{02}$ and $p_{02}.$ Substituting for $\alpha _{1}$
and $\alpha _{2}$ from (\ref{q2}) in the transformation (\ref{tr2}) and
using the method of undetermined coefficients, we find\newline
\newline
\indent $p_{02}=d_{1}\frac{d_{1}^{3}c_{02}-d_{1}^{2}\left(
c_{1}d_{02}+c_{1}l_{02}-l_{1}c_{02}\right) -d_{2}^{2}\left(
c_{1}d_{20}+c_{1}l_{20}-l_{1}c_{20}\right) +d_{1}d_{2}\left(
c_{1}d_{11}-d_{1}c_{11}+d_{2}c_{20}+c_{1}l_{11}-l_{1}c_{11}\right) }{%
c_{1}\left( d_{1}l_{2}-d_{2}l_{1}\right) ^{3}}$

\noindent and

$s_{02}=\allowbreak \frac{-\left( d_{2}c_{02}+l_{2}c_{02}\right)
d_{1}^{3}+\left( c_{1}d_{02}+d_{2}c_{11}+c_{1}l_{02}+l_{2}c_{11}\right)
\allowbreak d_{1}^{2}d_{2}-\left(
c_{1}d_{11}+d_{2}c_{20}+c_{1}l_{11}+l_{2}c_{20}\right) \allowbreak
d_{1}d_{2}^{2}+c_{1}d_{2}^{3}\left( d_{20}+l_{20}\right) }{c_{1}\left(
d_{1}l_{2}-d_{2}l_{1}\right) ^{3}}.$

These yield 
\begin{equation}
m_{2}=-\frac{%
c_{02}d_{1}^{3}-c_{11}d_{1}^{2}d_{2}-c_{1}d_{02}d_{1}^{2}+c_{20}d_{1}d_{2}^{2}+c_{1}d_{11}d_{1}d_{2}-c_{1}d_{20}d_{2}^{2}%
}{c_{1}\left( d_{1}l_{2}-d_{2}l_{1}\right) ^{2}}.  \label{m2}
\end{equation}%
Notice that, the both curves $\widehat{\beta }_{1,2}\left( \mu \right) =0$
are tangent to the $\mu _{2}-$axis at the origin and they are parabola-like
curves. $\square $

\begin{remark}\label{rem1} The new non-degeneracy condition (\ref{newt}) does not use any
coefficient from $\beta _{2}\left( \alpha \right)$ but only from $\beta _{1}\left( \alpha \right)$ and $L\left( \alpha \right).$
\end{remark}

The equation (\ref{y1y2}) becomes 
\begin{equation}
\widehat{L}_{2}\left( \mu \right) y^{2}+\widehat{\beta }_{2}\left( \mu
\right) y+\widehat{\beta }_{1}\left( \mu \right) =0,  \label{s2}
\end{equation}%
where $\widehat{L}_{2}\left( \mu \right) =L_{0}\left( 1+O\left( \left\vert
\mu \right\vert \right) \right) \neq 0.$ The new discriminant of (\ref{s2})
is $\widehat{\Delta }\left( \mu \right) =\Delta \circ S^{-1}\left( \mu
\right) ,$ which, by the transformation (\ref{tr2}), becomes 
\begin{equation}
\widehat{\Delta }\left( \mu \right) =\mu _{1},  \label{del}
\end{equation}%
while $\widehat{\beta }_{1,2}\left( \mu \right) $ are given in (\ref{be1})
and (\ref{be2}).

\begin{figure}[h]
\centering
\includegraphics[width=0.9\textwidth]{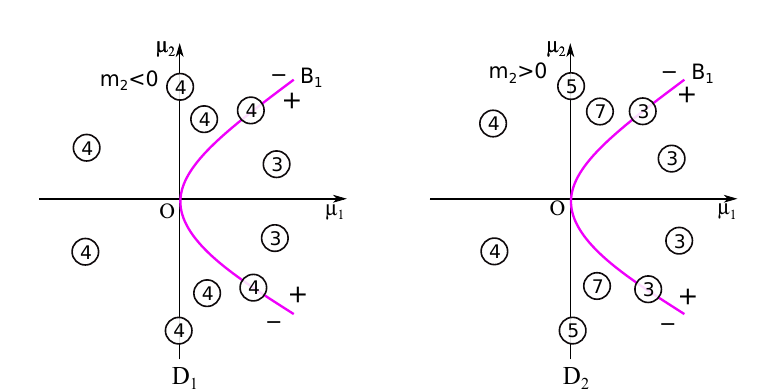}
\caption{Bifurcation diagrams corresponding to $L_0<0,$ $c_1d_1\neq 0,$
respectively, $m_2<0$ and $m_2>0.$}
\label{d1d2}
\end{figure}

\begin{remark}\label{rem2} 
The following theorem describes the bifurcation diagrams of the
system (\ref{tru})-(\ref{truf}) in the new parametric space $\mu _{1}O\mu
_{2}.$ In elaborating the diagrams from Figure \ref{d1d2} and Figure \ref%
{d3d4}, the curves $B_{1}$ and $B_{2}$ are approximated by $\left(
B_{1}\right) $ $\mu _{1}=m_{2}^{2}\mu _{2}^{4}$ and $\left( B_{2}\right) $ $%
\mu _{1}=4L_{0}\frac{c_{1}}{d_{1}}m_{2}\mu _{2}^{2},$ which are
parabola-like curves. Their relative positions one to another are given by $%
m_{2},$ $L_{0}$ and $\frac{c_{1}}{d_{1}}.$ By (\ref{del}), the curve $%
\widehat{\Delta }\left( \mu \right) =0$ coincides to the $\mu _{2}-$axes.
\end{remark}

\bigskip

\begin{figure}[h]
\centering
\includegraphics[width=0.9\textwidth]{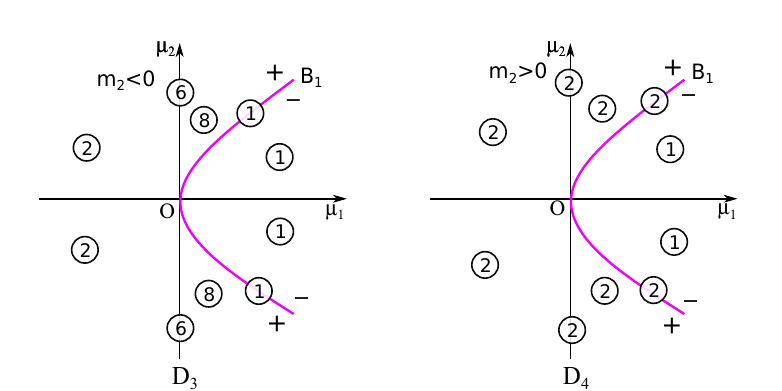}
\caption{Bifurcation diagrams corresponding to $L_0>0,$ $c_1d_1\neq 0,$
respectively, $m_2<0$ and $m_2>0.$}
\label{d3d4}
\end{figure}

\begin{theorem}\label{th3} 
The behavior of the system (\ref{tru})-(\ref{truf}) is described by
four bifurcation diagrams. More exactly, by Figure \ref{d1d2} if $L_0<0,$
respectively, Figure \ref{d3d4} if $L_0>0.$
\end{theorem}

\textit{Proof.} Let $\left( \mu _{1}^{B_{1}},\mu _{2}\right) \in B_{1}$ and $%
\left( \mu _{1}^{B_{2}},\mu _{2}\right) \in B_{2}$ be two points from the
two curves. Then, (\ref{b1}) and (\ref{b2}) yield 
\begin{equation}
\mu _{1}^{B_{2}}-\mu _{1}^{B_{1}}=k_{1}\mu _{2}^{2}\left( 1+O\left( \mu
_{2}\right) \right) \neq 0,  \label{comp}
\end{equation}%
for $\left\vert \mu _{2}\right\vert $ sufficiently small, $\mu _{2}\neq 0,$
where $k_{1}=4L_{0}\frac{c_{1}}{d_{1}}m_{2}.$

Notice that $B_{1}$ lies in $\left\{ \mu _{1}>0\right\} $ for all $m_{2}\neq
0,$ while $B_{2}\subset \left\{ \mu _{1}>0\right\} $ if $k_{1}>0,$
respectively, $B_{2}\subset \left\{ \mu _{1}<0\right\} $ if $k_{1}<0.$

Assume first $L_{0}<0,$ $m_{2}<0$ and $c_{1}d_{1}<0,$ thus, $k_{1}<0.$ Then $%
B_{2}\subset \left\{ \mu _{1}<0\right\} $ and $\widehat{\beta }_{1}\left(
\mu \right) =-\frac{1}{4L_{0}}\left( \mu _{1}-m_{2}^{2}\mu _{2}^{4}\right)
>0 $ inside the parabola $B_{1}$ (corresponding to $\mu _{1}>0$) and $%
\widehat{\beta }_{1}\left( \mu \right) <0$ on the exterior of $B_{1},$
Figure \ref{d1d2}. It follows from the Table 1 that, the phase portrait is $%
4 $ whenever $L_{0}<0$ and $\widehat{\beta }_{1}\left( \mu \right) <0,$ for
any signs of $\widehat{\beta }_{2}\left( \mu \right) $ and $\widehat{\Delta }%
\left( \mu \right) ,$ including on $\widehat{\Delta }\left( \mu \right) =0.$
Notice that $\widehat{\beta }_{2}\left( \mu \right) =-\frac{d_{1}}{%
4L_{0}c_{1}}\mu _{1}+m_{2}\mu _{2}^{2}<0$ on $\widehat{\Delta }\left( \mu
\right) =\mu _{1}>0.$ On the other hand, the phase portrait is $3$ if $%
L_{0}<0,$ $\widehat{\beta }_{1}\left( \mu \right) >0$ and $\widehat{\Delta }%
\left( \mu \right) =\mu _{1}>0,$ independent on the sign of $\widehat{\beta }%
_{2}\left( \mu \right) ,$ by Table 1. These prove the bifurcation diagram is 
$D_{1}$ from Figure \ref{d1d2}.

Assume now $c_{1}d_{1}>0,$ while $L_{0}<0$ and $m_{2}<0.$ Then $%
B_{1},B_{2}\subset \left\{ \mu _{1}>0\right\} $ and, by (\ref{comp}), $\mu
_{1}^{B_{2}}>\mu _{1}^{B_{1}}>0,$ thus, the parabola $B_{2}$ lies in the
interior of the parabola $B_{1},$ where $\widehat{\beta }_{1}\left( \mu
\right) >0$ and $\widehat{\Delta }\left( \mu \right) =\mu _{1}>0.$ Thus, by
Table 1, the phase portrait is $3$ in the interior of $B_{1}$ for any sign
of $\widehat{\beta }_{2}\left( \mu \right) ,$ respectively, $4$ in the
exterior of $B_{1}$ because $\widehat{\beta }_{2}\left( \mu \right) =-\frac{%
d_{1}}{4L_{0}c_{1}}\mu _{1}+m_{2}\mu _{2}^{2}<0$ on the exterior of $B_{1}.$
Thus, the same bifurcation diagram $D_{1},$ Figure \ref{d1d2}, characterizes
the dynamics of the system (\ref{tru})-(\ref{truf}) in this case.

In the second case, assume $L_{0}<0$ and $m_{2}>0,$ and consider first $%
c_{1}d_{1}>0,$ thus, $k_{1}<0.$ Then $B_{2}\subset \left\{ \mu
_{1}<0\right\} $ and $\widehat{\beta }_{1}\left( \mu \right) >0$ inside of $%
B_{1}$ ( $\mu _{1}>0$ ) and $\widehat{\beta }_{1}\left( \mu \right) <0$ on
the exterior of $B_{1},$ Figure \ref{d1d2}$\left( D_{2}\right) .$ By Table
1, the phase portrait is $4$ whenever $L_{0}<0,$ $\widehat{\Delta }\left(
\mu \right) =\mu _{1}<0$ and $\widehat{\beta }_{1}\left( \mu \right) <0,$
for any signs of $\widehat{\beta }_{2}\left( \mu \right) .$ On $\widehat{%
\Delta }\left( \mu \right) =0,$ the phase portrait becomes $5$ since $%
\widehat{\beta }_{2}\left( \mu \right) >0$ on $\mu _{1}\geq 0,$ while on the
exterior of $B_{1}$ and $\mu _{1}>0,$ it is $7$ because $\widehat{\beta }%
_{1}\left( \mu \right) <0.$ On $B_{1}$ and the interior of $B_{1},$ it
becomes $3,$ because $\widehat{\beta }_{1,2}\left( \mu \right) >0$ and $%
\widehat{\Delta }\left( \mu \right) >0.$ Thus, the bifurcation diagram
corresponding to this case is $D_{2}$ from Figure \ref{d1d2}.

When $c_{1}d_{1}<0,$ $L_{0}<0$ and $m_{2}>0,$ the curves $B_{1,2}\subset
\left\{ \mu _{1}>0\right\} $ with $k_{1}>0.$ By (\ref{comp}), $\mu
_{1}^{B_{2}}>\mu _{1}^{B_{1}}>0,$ thus, the parabola $B_{2}$ lies again in
the interior of the parabola $B_{1},$ where $\widehat{\beta }_{1}\left( \mu
\right) >0$ and $\widehat{\Delta }\left( \mu \right) >0.$ It follows from
Table 1 that, the phase portrait is $3$ inside of $B_{1},$ independent on
the sign of $\widehat{\beta }_{2}\left( \mu \right) ,$ while it remains $3$
on $B_{1}$ because $\widehat{\beta }_{2}\left( \mu \right) >0$ on the
exterior of $B_{2}.$ Further, on the exterior of $B_{1}$ and $\mu _{1}>0,$
the phase portrait becomes $7,$ since $L_{0}<0,$ $\widehat{\beta }_{1}\left(
\mu \right) <0,$ $\widehat{\Delta }\left( \mu \right) >0$ and $\widehat{%
\beta }_{2}\left( \mu \right) >0,$ which transforms in $5$ on $\mu _{1}=0,$
respectively, in $4$ on $\mu _{1}<0,$ because $L_{0}<0,$ $\widehat{\beta }%
_{1}\left( \mu \right) <0,$ $\widehat{\Delta }\left( \mu \right) <0$ and $%
\widehat{\beta }_{2}\left( \mu \right) >0.$ Therefore, the bifurcation
diagram corresponding to this case is $D_{2}$ as well.

One can proceed similarly for the case $L_{0}>0.$ Two bifurcation diagrams, $%
D_{3}$ and $D_{4}$ from Figure \ref{d3d4}, describe the dynamics of the
system (\ref{tru})-(\ref{truf}) in this case. $\square $

\begin{figure}[h!]
\begin{center}
\includegraphics[width=0.45\textwidth]{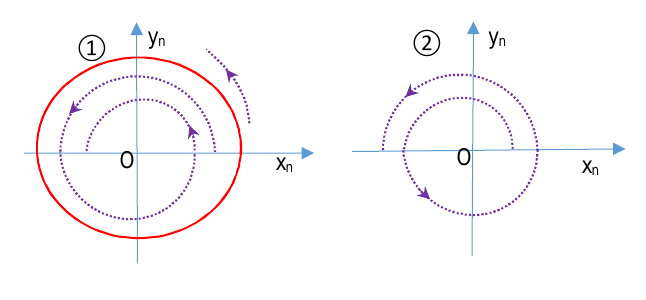} \includegraphics[width=0.45%
\textwidth]{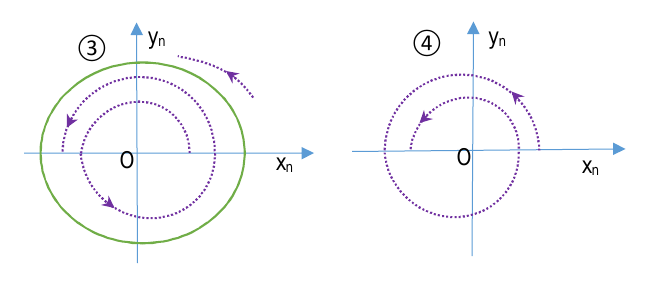} \includegraphics[width=0.45\textwidth]{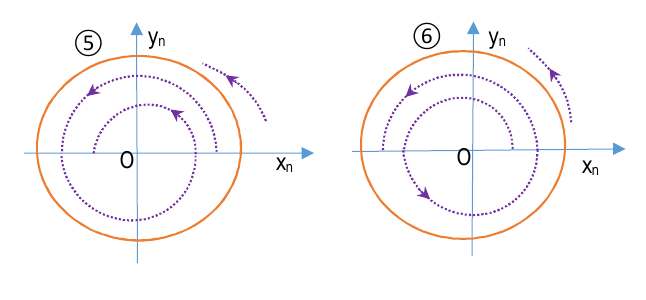} %
\includegraphics[width=0.45\textwidth]{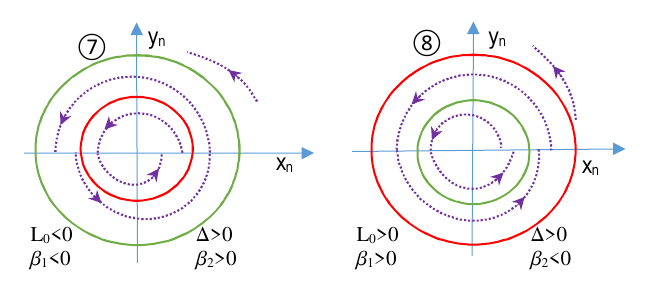}
\end{center}
\caption{Generic phase portraits of the system (\protect\ref{tru})-(\protect
\ref{truf}).}
\label{f1}
\end{figure}

The bifurcation regions we encounter in our study are among the ones given
in Table 1, which were described in Theorem \ref{th3}. Figure \ref{f1}
illustrates the possible generic phase portraits corresponding to the
regions from the two tables.

\newpage

\begin{center}
Table 1. \textit{The regions in the parametric plane $\mu _{1}\mu _{2}$
defined by $\widehat{\Delta }(\mu ),$ $\widehat{\beta }_{1,2}(\mu )$ and $%
L_{0}.$}

\begin{tabular}{lllll|llllll}
$L_{0}$ & $\widehat{\Delta }\left( \mu \right) $ & $\widehat{\beta }_{1}$ & $%
\widehat{\beta }_{2}$ & $Region$ &  & $L_{0}$ & $\widehat{\Delta }\left( \mu
\right) $ & $\widehat{\beta }_{1}$ & $\widehat{\beta }_{2}$ & $Region$ \\ 
\hline
$-$ & $-$ & $-$ & $\pm ,0$ & $4$ &  & $+$ & $+$ & $-$ & $\pm ,0$ & $1$ \\ 
$-$ & $+$ & $-$ & $-$ & $4$ &  & $+$ & $+$ & $0$ & $-$ & $1$ \\ 
$-$ & $0$ & $-$ & $-$ & $4$ &  & $+$ & $-$ & $+$ & $\pm ,0$ & $2$ \\ 
$-$ & $+$ & $0$ & $-$ & $4$ &  & $+$ & $+$ & $+$ & $+$ & $2$ \\ 
$-$ & $0$ & $0$ & $0$ & $4$ &  & $+$ & $0$ & $+$ & $+$ & $2$ \\ 
$-$ & $+$ & $+$ & $\pm ,0$ & $3$ &  & $+$ & $+$ & $0$ & $+$ & $2$ \\ 
$-$ & $+$ & $0$ & $+$ & $3$ &  & $+$ & $0$ & $0$ & $0$ & $2$ \\ 
$-$ & $+$ & $-$ & $+$ & $7$ &  & $+$ & $0$ & $+$ & $-$ & $6$ \\ 
$-$ & $0$ & $-$ & $+$ & $5$ &  & $+$ & $+$ & $+$ & $-$ & $8$%
\end{tabular}
\end{center}

\begin{example} Consider a two-dimensional map $\left( \rho _{n},\varphi _{n}\right)
\longmapsto \left( \rho _{n+1},\varphi _{n+1}\right) $ given in polar
coordinates by 
\begin{equation}
\left\{ 
\begin{array}{cc}
\rho _{n+1}= & \rho _{n}\left( 1+\beta _{1}\left( \alpha \right) +\beta
_{2}\left( \alpha \right) \rho _{n}^{2}+\allowbreak L_{2}\left( \alpha
\right) \rho _{n}^{4}\right) \\ 
\\
\varphi _{n+1}= & \varphi _{n}+\theta _{0}\ \ \ \  \ \ \ \ \ \ \ \ \ \ \ \ \ \ \ \ \ \ \ \ \ \ \ \ \ \ \ \  \ \ \ \  \ \text{ }%
\end{array}%
\right. ,  \label{ex1}
\end{equation}
where $\theta _{0}$ is fixed, $0<\theta _{0}\ <\pi ,$ $\beta _{1}\left(
\alpha \right) =\alpha _{1}+\alpha _{2}+2\alpha _{1}^{2}+\alpha _{2}^{2},$ $%
\beta _{2}\left( \alpha \right) =\alpha _{1}+\alpha _{2}+2\alpha _{1}\alpha
_{2}$ and $L_{2}\left( \alpha \right) =1+\alpha _{1}+2\alpha _{2}+\alpha
_{1}^{2}+\alpha _{2}^{3}.$
\end{example}

\begin{remark} The map (\ref{ex1}) is degenerate with respect to the change of parameters $%
\left( \alpha _{1},\alpha _{2}\right) \longmapsto \left( \beta _{1},\beta
_{2}\right) $ since $c_{1}d_{2}-c_{2}d_{1}=0,$ thus, it cannot be studied with the known methods for the non-degenerate case. On the other hand, the
transformation (\ref{tr2}) proposed in this paper, $\left( \alpha _{1},\alpha _{2}\right) \overset{S}{%
\longmapsto }\left( \mu _{1},\mu _{2}\right),$ is
regular at $\left( 0,0\right) $ since (\ref{newt}) is satisfied, thus, it can be applied to study the degenerate Chenciner bifurcation in this map.
\end{remark}

The inverse transformation $S^{-1}$ in its linear terms is given by $\alpha
_{1}=-\frac{3}{4}\mu _{1}-\mu _{2}$ and $\alpha _{2}=\frac{1}{2}\mu _{1}+\mu
_{2}.$ In order to find $m_{2}$ without applying (\ref{m2}), we need $S^{-1}$
in its linear and quadratic terms, which is of the form 
\begin{equation*}
\alpha _{1}=-\frac{3}{4}\mu _{1}-\mu _{2}+s_{20}\mu _{1}^{2}+s_{11}\mu
_{1}\mu _{2}+s_{02}\mu _{2}^{2}\ \ \ \ \text{ and }\ \ \ \ \alpha _{2}=\frac{%
1}{2}\mu _{1}+\mu _{2}+p_{20}\mu _{1}^{2}+p_{11}\mu _{1}\mu _{2}+p_{02}\mu
_{2}^{2}.
\end{equation*}

Using the method of undetermined coefficients in (\ref{tr2}), we obtain a
linear system in the unknowns $s_{ij}$ and $p_{ij},$ $i+j=2,$ whose solution
is $p_{02}=7,$ $p_{11}=\frac{17}{2},$ $p_{20}=\frac{89}{32},$ $s_{02}=-10,$ $%
s_{11}=-\frac{49}{4}$ and $s_{20}=-\frac{261}{64}\allowbreak .$ These lead
to $\widehat{L}_{2}\left( \mu \right) =1+\allowbreak \frac{1}{4}\mu _{1}+\mu
_{2}+5\mu _{2}^{2}$ and $\widehat{\beta }_{2}\left( \mu \right) =-\frac{1}{4}%
\mu _{1}-5\allowbreak \mu _{2}^{2}$ in their lowest terms, thus, $m_{2}=-5.$

We notice that $m_{2}=-5$ could be obtained directly from the coefficients
of $\beta _{1,2}\left( \alpha \right) $ and $L_{2}\left( \alpha \right) $ by
formula (\ref{m2}).

Using $S^{-1}$ up to quadratic terms is not sufficient for finding $\widehat{%
\beta }_{1}\left( \mu \right) .$ However, $\widehat{\beta }_{1}\left( \mu
\right) $ can be determined from the substitution (\ref{tr2}) by $\widehat{%
\beta }_{2}^{2}\left( \mu \right) -4\widehat{\beta }_{1}\left( \mu \right) 
\widehat{L}_{2}\left( \mu \right) -\mu _{1}=0,$ which yields $\widehat{\beta 
}_{1}\left( \mu \right) =-\frac{1}{4}\mu _{1}+25\allowbreak \mu _{2}^{4}$ in
their lowest terms. Since $L_{0}=1$ and $m_{2}=-5,$ the corresponding
bifurcation diagram of the map (\ref{ex1}) is $D_{3},$ Figure \ref{d3d4}.

Let us illustrate numerically the behavior of the map (\ref{ex1}) for
different values of $\alpha =\left( \alpha _{1},\alpha _{2}\right) $
corresponding to the four different regions of $D_{3},$ which we denote by $%
R_{1},$ $R_{2},$ $R_{6}$ and $R_{8}.$

To this end, we proceed as it follows. For a given numerical value $\left(
\alpha _{1},\alpha _{2}\right) ,$ we determine $\mu =\left( \mu _{1},\mu
_{2}\right) $ by the transformation $S$ from $\mu _{1}=\beta _{2}^{2}\left(
\alpha \right) -4\beta _{1}\left( \alpha \right) L_{2}\left( \alpha \right) $
and $\mu _{2}=\beta _{2}\left( \alpha \right) +L_{2}\left( \alpha \right)
-L_{0},$ and find the region of $D_{3}$ where $\mu $ lies. Thus, the
behavior of the map (\ref{ex1}) should be in agreement with the
corresponding region. To probe this, we integrate numerically the map in
Matlab, with $\left( \alpha _{1},\alpha _{2}\right) $ fixed in the begining,
and find different orbits $\left( x_{n},y_{n}\right) ,$ where $x_{n}=\rho
_{n}\cos \varphi _{n}$ and $y_{n}=\rho _{n}\sin \varphi _{n},$ for $n$
taking all integer values from $1$ to a fixed value $N.$ The cartesian
values $x_{n}$ and $y_{n}$ are then plotted in the same diagram.

\begin{figure}[h]
\begin{center}
\includegraphics[width=0.4\textwidth]{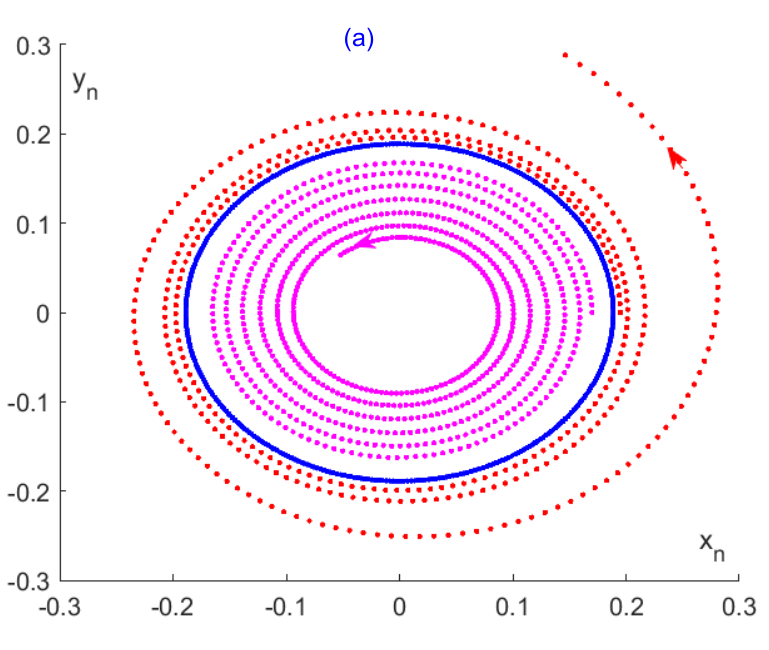} \includegraphics[width=0.4%
\textwidth]{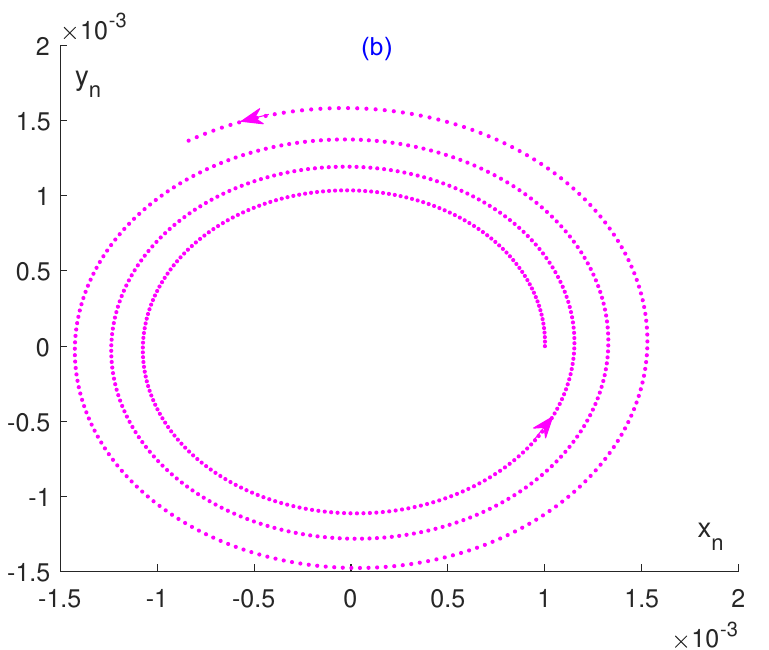} \includegraphics[width=0.4\textwidth]{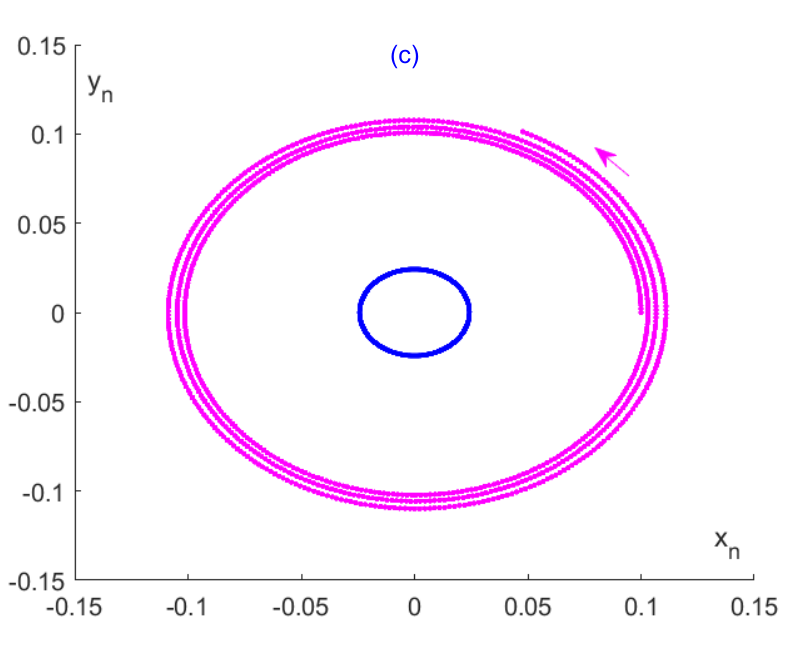} %
\includegraphics[width=0.4\textwidth]{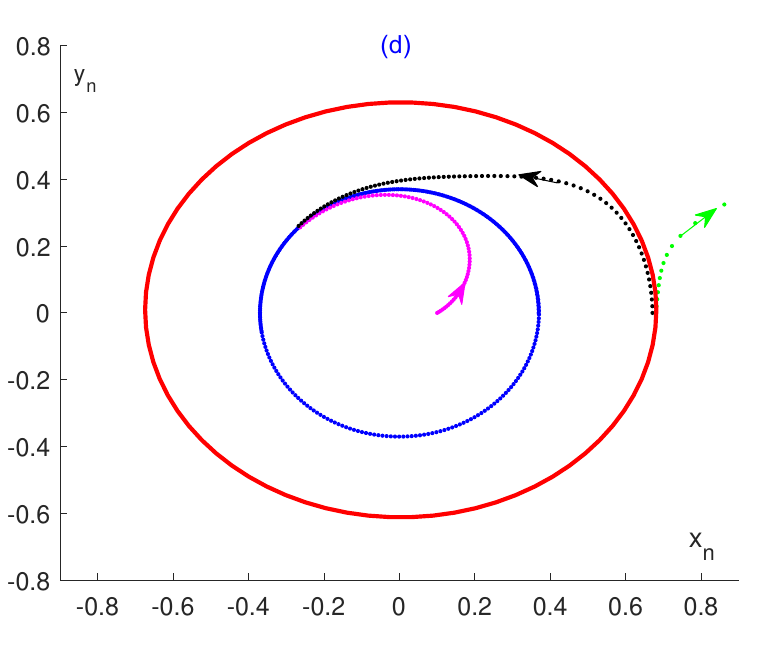}
\end{center}
\caption{The behavior of the map (\protect\ref{ex1}) on the four regions of
the bifurcation diagram $D_3.$}
\label{f2}
\end{figure}

Consider first $\alpha _{1}=-0.017$ and $\alpha _{2}=0.015,$ which lead to $%
\mu _{1}=\Delta \left( \alpha \right) =\allowbreak 4.8579\times 10^{-3}$ and 
$\mu _{2}=\allowbreak 1.0782\times 10^{-2},$ respectively, $\widehat{\beta }%
_{1}\left( \mu \right) =\allowbreak -1.\allowbreak 2141\times 10^{-3}$ and $%
\widehat{\beta }_{2}\left( \mu \right) =-1.7958\times 10^{-3},$ thus, $%
\left( \mu _{1},\mu _{2}\right) \in R_{1}\in D_{3},$ with $\sqrt{y_{1}}%
=\allowbreak 0.18876.$ For these values and $\theta _{0}=0.05,$ three orbits
are obtained and presented in Figure \ref{f2}(a). The first orbit (magenta)
starts at $\left( \rho _{1},\varphi _{1}\right) =\left( 0.17,0\right) $ and
tends to the origin $\left( 0,0\right) $ as $n$ increases from $1$ to $%
N=800. $ The second orbit (blue) approximates an invariant closed curve
(which is a circle); it starts at $\left( \rho _{1},\varphi _{1}\right)
=\left( 0.18876,0\right) $ and was obtained with $N=400.$ The third orbit
(red) starts at $\left( 0.195,0\right) ;$ it departs from the invariant
closed curve for $n$ increasing and may escape to infinity. It follows that,
the closed invariant circle is unstable. We notice that, the three orbits
from Figure \ref{f2}(a) are in perfect agreement with the phase portrait $1$
presented schematically in Figure \ref{f1}.

For the second region $R_{2},$ let $\alpha _{1}=-0.015$ and $\alpha
_{2}=0.015,$ which lead to $\mu _{1}=\Delta \left( \alpha \right)
=\allowbreak \allowbreak -2.7\times 10^{-3}$ and $\mu _{2}=\allowbreak
1.\,\allowbreak 4778\times 10^{-2},$ respectively, $\widehat{\beta }%
_{1}\left( \mu \right) =\allowbreak \allowbreak 6.8\times 10^{-4}$ and $%
\widehat{\beta }_{2}\left( \mu \right) =-4\times 10^{-4}<0,$ thus, $\left(
\mu _{1},\mu _{2}\right) \in R_{2}\in D_{3}.$ An orbit for these values and $%
\theta _{0}=0.03,$ $N=700,$ $\left( \rho _{1},\varphi _{1}\right) =\left(
0.001,0\right) ,$ is presented in Figure \ref{f2}(b). The orbit departs from
the origin and may escape to infinity, which is in agreement with the phase
portrait $2$ from Figure \ref{f1}. We notice that $R_{2}$ may contain points 
$\left( \mu _{1},\mu _{2}\right) $ with $\widehat{\beta }_{2}\left( \mu
\right) >0,$ which occurs, for example, at $\alpha _{1}=0.015$ and $\alpha
_{2}=0.015,$ which yield $\mu _{1}=-0.127$ and $\mu _{2}=0.0\allowbreak 75,$
respectively, $\widehat{\beta }_{1}\left( \mu \right) =\allowbreak
\allowbreak 0.032$ and $\widehat{\beta }_{2}\left( \mu \right) =0.0031.$ The
orbits simulated for these values are also in agreement with Figure \ref{f1}.

In the third case, let $\alpha _{1}=-0.015719$ and $\alpha _{2}=0.015,$
which lead to $\mu _{1}=\Delta \left( \alpha \right) =\allowbreak
\allowbreak 0$ and $\mu _{2}=1.3341\times 10^{-2},$ respectively, $\widehat{%
\beta }_{1}\left( \mu \right) =\allowbreak 6.1\times 10^{-7},$ $\widehat{%
\beta }_{2}\left( \mu \right) =-8.9\times 10^{-4}$ and $\sqrt{y_{1,2}}%
=\allowbreak 0.0242,$ thus, $\left( \mu _{1},\mu _{2}\right) \in R_{6}\in
D_{3}.$ Setting $\theta _{0}=0.02$ and $N=1000,$ two orbits are depicted in
Figure \ref{f2}(c). One (blue) starts at $\left( \rho _{1},\varphi
_{1}\right) =\left( 0.024223,0\right) $ and approximates a closed invariant
orbit (circle), while the other (magenta) starts at $\left( 0.1,0\right) $
and departs from the closed invariant orbit. Orbits starting from the
interior of the invariant circle tend slowly to the circle for $n$
increasing. Thus, the invariant circle is stable from the interior and
unstable from the exterior. This behavior is in agreement to the phase
portrait $6$ from Figure \ref{f1}.

Finally, let $\alpha _{1}=-0.5$ and $\alpha _{2}=0.05.$ Then $\Delta \left(
\alpha \right) >0,$ $\rho _{1}=\sqrt{y_{1}}=\allowbreak 0.6718$ and $\rho
_{2}=\sqrt{y_{2}}=0.3699,$ which, by $S,$ lead to $\mu _{1}=\Delta \left(
\alpha \right) =\allowbreak 0.0714$ and $\mu _{2}=\allowbreak -0.6498,$
respectively, $\widehat{\beta }_{1}\left( \mu \right) =\allowbreak 4.4$ and $%
\widehat{\beta }_{2}\left( \mu \right) \allowbreak =-2.1,$ that is, $\left(
\mu _{1},\mu _{2}\right) \in R_{8}$ from $D_{3}.$ Numerical simulations of
this case with $\theta _{0}=0.03$ are presented in Figure \ref{f2}(d).
Notice the appearance of two invariant closed curves, one stable (blue) and
the other (red) unstable. This is in perfect agreement with the phase
portrait $8,$ presented schematically in Figure \ref{f1}.

\section{Conclusions}

The study we performed in this work brings to light new generic properties
of a degenerate form of the Chenciner bifurcation. The degeneracy we tackled
refers at the regularity of the transformation of parameters, which is used
for obtaining a normal form in nondegenerate case. We studied the
bifurcation when the transformation of parameters is not regular at $(0,0),$
thus, the classical results can not be applied. We proposed a different
change of parameters, which was proved to be successful for approaching the
degeneracy. More exactly, when the degeneracy condition occurs, the new
transformation is regular at $(0,0)$ and, more importantly, can be used for
studying the bifurcation. The cost of using this transformation is the
appearance of a different generic condition. The new condition uses
different coefficients to exist than the classical one, thus, they are
complementary one to another. We exemplified the application of the method
developed in this work on a particular map, which cannot be studied with the
method used for nondegenerate case. When the both generic conditions fail, a
new study is needed for exploring the behavior of the Chenciner bifurcation.
This remains an open problem.

\section{Acknowledgments}

This research was supported by Horizon2020-2017-RISE-777911 project and a
grant of the Ministry of Education and Research, CNCS/CCCDI - UEFISCDI,
PN-III-P3-3.6-H2020-2020-0100.

\end{document}